\documentclass[11pt]{article}
\usepackage{amscd}
\usepackage{amsfonts}
\usepackage{amsmath}
\usepackage{amssymb}
\usepackage{amsthm}
\usepackage{bbm}
\usepackage{CJK}
\usepackage{fancyhdr}
\usepackage{graphicx}
\usepackage{hyperref}
\usepackage{indentfirst}
\usepackage{latexsym}
\usepackage{mathrsfs}
\usepackage{xypic}

\usepackage[top=1in,bottom=1in,left=1.25in,right=1.25in]{geometry}
\def\<{\langle}
\def\>{\rangle}
\def\a{\alpha}
\def\b{\beta}

\def\c{\cdot}

\def\D{\Delta}

\def\g{\gamma}

\def\o{\otimes}

\def\v{\varepsilon}

\begin{document}
\renewcommand{\baselinestretch}{1.2}
\renewcommand{\arraystretch}{1.0}
\title{\textbf{A Maschke type theorem for relative Hom-Hopf modules  }}
\author{{Shuang-jian Guo$^{1}$ \quad   Xiu-li Chen$^{2}$\thanks{Corresponding author, email: xiulichen1021@126.com }}\\
{\small 1. School of Mathematics and Statistics, Guizhou University of Finance and Economics}\\
{\small Guiyang,  550025, P. R. China}\\
{\small 2. Department of Mathematics, Southeast  University}\\
 {\small Nanjing, 210096,  P. R. China}}
 \date{}
 \maketitle
\noindent\textbf {\bf Abstract} In this paper, we prove a Maschke type theorem for the category of relative Hom-Hopf modules. In fact,  we give necessary and sufficient conditions for the functor that forgets the $(H, \a)$-coaction to be separable. This leads to a generalized notion of integrals.


\noindent\textbf {MR(2010)Subject Classification} 16T05

\section{Introduction}
The present paper investigates variations on the theme of
Hom-algebras, a topic which has recently received much attention
from various researchers. The study of Hom-associative algebras
originates with work by Hartwig, Larsson and Silvestrov in the Lie
case [9], where a notion of Hom-Lie algebra was introduced in the
context of studying deformations of Witt and Virasoro algebras.
Later, it was extended to the associative case by Makhlouf and
Silverstrov in [10]. Now the associativity is replaced by
Hom-associativity $\alpha(a)(bc)=(ab)\alpha(c)$. Hom-coassociativity
for a Hom-coalgebra can be considered in a similar way, see [10].
Caenepeel etc.  [1] studied Hom-structures from the point of view of
monoidal categories. This leads to the natural definition of monoidal Hom-algebras,
Hom-coalgebras, etc. They  constructed a symmetric monoidal category, and then
introduced monoidal Hom-algebras, Hom-coalgebras, etc. as algebras, coalgebras,
etc. in this monoidal category.

The notion of a relative $(H,B)$-Hopf module, where $H$ is a Hopf algebra over a field $k$ and $B$
is a right coideal subalgebra of $H$, was introduced and studied by Takeuchi in [12]. Later,
in [3](see also [4]),   Doi noted that the notion of an $(H,B)$-Hopf module works well if $B$ is a
right $H$-comodule algebra, Using this module, he proved that the existence of a total integral
$\phi: H\rightarrow B$ is equivalent to $B$ being a relative injective $H$-comodule, and it is also equivalent to
any$(H,B)$-Hopf module $M$ being a relative injective $H$-comodule in [5]. Also, in [5], using a
commutative assumption for $H$, he deduced a version of the Maschke type theorem for $(H,B)$-Hopf
modules which states that every exact sequence of $(H,B)$-Hopf modules which splits $B$-linearly,
also splits $(H,B)$-linearly. Afterwards, Doi proved in [5] that the commutative condition can
be removed and replaced by some technical conditions involving the center of $B$.  Caenepeel etc. [2] proved a Maschke type theorem for the category of relative Hopf modules. In
fact, they gave necessary and sufficient conditions for the functor that forgets the $H$-coaction to be separable. This leads to a generalized notion of integrals of Doi [5].

In this paper we study the generalization of the previous results to the Hom-Hopf algebras.
In Sec.3, we introduce the notion of a relative Hom-Hopf module and prove that the functor $F $ from  the category of relative Hom-Hopf modules to the category of right $(A, \b)$-Hom-modules has a right
 adjoint(see Proposition 3.3).
In Sec.4,  we introduce the notion of total integrals for
Hom-comodule algebras, which is an effective tool for
investigating properties of relative Hom-Hopf modules. As an
important application, we investigate the injectivity of relative
Hom-Hopf modules(see Proposition 4.3), which generalizes the main
result in [3]. In Sec.5,  we obtain the main result of this paper.
 we give necessary and sufficient conditions for the functor that forgets the $(H,\a)$-coaction to be separable (see Theorem 5.2), we prove a Maschke type theorem for the category of relative Hom-Hopf modules as an application. In fact, let $(A,\b)$ be a right $(H,\a)$-Hom-comodule algebra
with a total integral $\phi: (H,\a)\rightarrow (A,\b)$. If $\phi: (H,\a)\rightarrow
(Z(A),\b)$ (the center of $A$) is a multiplication map, then every
short exact sequence of relative Hom-Hopf modules
$$0\longrightarrow (M,\mu)\stackrel{f}\longrightarrow
(N,\nu)\stackrel{g}\longrightarrow (P,\pi)\longrightarrow 0$$ which splits as a
sequence of $(A,\b)$-Hom-modules also splits as a sequence of
relative Hom-Hopf modules.

\section{Preliminaries}
Throughout this paper we work over a commutative ring $k$, we recall from [1]
some information about Hom-structures which are needed in what follows.

Let $\mathcal{C}$ be a category. We introduce a new category
$\widetilde{\mathscr{H}}(\mathcal{C})$ as follows: objects are
couples $(M, \mu)$, with $M \in \mathcal{C}$ and $\mu \in
Aut_{\mathcal{C}}(M)$. A morphism $f: (M, \mu)\rightarrow (N, \nu )$
is a morphism $f : M\rightarrow N$ in $\mathcal{C}$ such that $\nu
\circ f= f \circ \mu$.

Let $\mathscr{M}_k$ denotes the category of $k$-modules.
~$\mathscr{H}(\mathscr{M}_k)$ will be called the Hom-category
associated to $\mathscr{M}_k$. If $(M,\mu) \in \mathscr{M}_k$, then
$\mu: M\rightarrow M$ is obviously a morphism in
~$\mathscr{H}(\mathscr{M}_k)$. It is easy to show that
~$\widetilde{\mathscr{H}}(\mathscr{M}_k)$ =
(~$\mathscr{H}(\mathscr{M}_k),~\otimes,~(I, I),~\widetilde{a},
~\widetilde{l},~\widetilde{r}))$ is a monoidal category by
Proposition 1.1 in [1]: the tensor product of $(M,\mu)$ and $(N,
\nu)$ in ~$\widetilde{\mathscr{H}}(\mathscr{M}_k)$ is given by the
formula $(M, \mu)\otimes (N, \nu) = (M\otimes N, \mu \otimes \nu)$.

Assume that $(M, \mu), (N, \nu), (P,\pi)\in
\widetilde{\mathscr{H}}(\mathscr{M}_k)$. The associativity and unit
constraints are given by the formulas
\begin{eqnarray*}
\widetilde{a}_{M,N,P}((m\o n)\o p)=\mu(m)\o (n\o \pi^{-1}(p)),\\
\widetilde{l}_{M}(x\o m)=\widetilde{r}_{M}(m\o x)=x\mu(m).
\end{eqnarray*}
An algebra in $\widetilde{\mathscr{H}}(\mathscr{M}_k)$ will be called a monoidal Hom-algebra.\\

\noindent{\bf Definition 2.1.} A monoidal Hom-algebra  is an
object is a triple $(A, \alpha)\in
\widetilde{\mathscr{H}}(\mathscr{M}_k)$ together with a $k$-linear
map $m_A: A\o A\rightarrow A$ and an element $1_A\in A$ such that
\begin{eqnarray*}
\alpha(ab)=\alpha(a)\alpha(b); ~\a(1_A)=1_A,\\
\alpha(a)(bc)=(ab)\alpha(c); ~a1_A=1_Aa=\a(a),
\end{eqnarray*}
for all $a,b,c\in A$. Here we use the notation $m_A(a\o b)=ab$.

\noindent{\bf Definition 2.2.} A monoidal Hom-coalgebra  is
an object $(C,\g)\in \widetilde{\mathscr{H}}(\mathscr{M}_k)$
together with $k$-linear maps $\Delta:C\rightarrow C\otimes C,~~~
\D(c)=c_{(1)}\o c_{(2)}$ (summation implicitly understood) and
$\gamma:C\rightarrow C$ such that
\begin{eqnarray*}
\D(\g(c))=\g(c_{(1)})\o \g(c_{(2)}); ~~~\varepsilon(\g(c))=\varepsilon(c),
\end{eqnarray*}
and
\begin{eqnarray*}
\g^{-1}(c_{(1)})\otimes c_{(2)(1)}\otimes c_{(2)(2)}=c_{(1)(1)}\otimes c_{(1)(2)}\otimes\g(c_{(2)}),~\varepsilon(c_{(1)})c_{(2)}=\varepsilon(c_{(2)})c_{(1)}=\g^{-1}(c)
\end{eqnarray*}
for all $c\in C$.\\

\noindent{\bf Definition 2.3.} A monoidal Hom-bialgebra
$H=(H,\a,m,\eta, \Delta,\varepsilon)$  is a bialgebra in
the symmetric monoidal category
$\widetilde{\mathscr{H}}(\mathscr{M}_k)$. This means that $(H, \a,
m,\eta)$ is a Hom-algebra, $(H,\Delta,\alpha)$ is a Hom-coalgebra
and that $\D$ and $\v$ are morphisms of Hom-algebras that is,
\begin{eqnarray*}
\Delta(ab)=a_{(1)}b_{(1)}\otimes a_{(2)}b_{(2)}; ~~~\Delta(1_H)=1_H\otimes 1_H,\\
\varepsilon(ab)=\varepsilon(a)\varepsilon(b),~\varepsilon(1_H)=1_H.
\end{eqnarray*}

\noindent{\bf Definition 2.4.} A monoidal Hom-Hopf algebra is
a monoidal Hom-bialgebra $(H, \alpha)$ together with a linear map
$S:H\rightarrow H$ in $\widetilde{\mathscr{H}}(\mathscr{M}_k)$ such
that
$$S\ast I=I\ast S=\eta\varepsilon,~S\alpha=\alpha S.$$

\noindent{\bf Definition 2.5.} Let $(A,\alpha)$ be a monoidal
Hom-algebra.  A right $(A,\alpha)$-Hom-module is an object
$(M,\mu)\in \widetilde{\mathscr{H}}(\mathscr{M}_k) $ consists of a
$k$-module and a linear map $\mu:M\rightarrow M$ together with a
morphism $\psi:M\otimes A\rightarrow M, \psi(m\c a)=m\c a$, in
$\widetilde{\mathscr{H}}(\mathscr{M}_k) $  such that
\begin{eqnarray*}
(m\c a)\c\alpha(b)=\mu(m)\c (ab);~~~~ m\c 1_A=\mu(m),
\end{eqnarray*}
for all $a\in A$ and $m\in M$. The fact that $\psi\in \widetilde{\mathscr{H}}(\mathscr{M}_k)$ means that
\begin{eqnarray*}
\mu(m\c a)=\mu(m)\c\alpha(a).
\end{eqnarray*}
A morphism $f: (M, \mu)\rightarrow (N, \nu)$ in $\widetilde{\mathscr{H}}(\mathscr{M}_k)$ is called right $A$-linear if it preserves the $A$-action, that is, $f(m\c a)=f(m)\c a$. $\widetilde{\mathscr{H}}(\mathscr{M}_k)_A$ will denote the category of right $(A, \a)$-Hom-modules and $A$-linear morphisms.\\

\noindent{\bf Definition 2.6.} Let $(C,\g)$ be a monoidal
Hom-coalgebra. A right $(C,\g)$-Hom-comodule  is an object
$(M,\mu)\in \widetilde{\mathscr{H}}(\mathscr{M}_k)$  together with a
$k$-linear map $\rho_M: M\rightarrow M\otimes C$ notation
$\rho_M(m)=m_{[0]}\o m_{[1]}$ in
$\widetilde{\mathscr{H}}(\mathscr{M}_k)$ such that
\begin{eqnarray*}
m_{[0][0]}\otimes (m_{[0][1]}\otimes
\g^{-1}(m_{[1]}))=\mu^{-1}(m_{[0]})\otimes \D_C(m_{[1]});
~m_{[0]}\varepsilon(m_{[1]})=\mu^{-1}(m),
\end{eqnarray*}
for all $m\in M$.  The fact that $\rho_M\in \widetilde{\mathscr{H}}(\mathscr{M}_k)$ means that
\begin{eqnarray*}
\rho_M(\mu(m))=\mu(m_{[0]})\otimes \g(m_{[1]}).
\end{eqnarray*}
Morphisms of right $(C, \g)$-Hom-comodule are defined in the obvious way. The
category of right $(C, \g)$-Hom-comodules will be denoted by $\widetilde{\mathscr{H}}(\mathscr{M}_k)^C$.

\section{Adjoint functor}

\noindent{\bf Definition 3.1.}  Let $(H,\alpha)$  be a monoidal Hom-Hopf algebra.  A monoidal Hom-algebra  $(A,\beta)$ is called  a right
$(H,\alpha)$-Hom-comodule algebra if $(A,\beta)$ is a right $(H,\alpha)$ Hom-comodule with coaction $\rho_A: A\rightarrow A\o H,
 ~~\rho_A(a)= a_{[0]}\o a_{[1]}$
 such that the following conditions satisfy,
 \begin{eqnarray*}
 &&\rho_{A} ( ab)= a_{[0]} b_{[0]} \otimes a_{[1]}b_{[1]},\\
 &&\rho_{A} ( 1_A)=1_A\o 1_H.
 \end{eqnarray*}
for all $a,b \in A$.\\

\noindent{\bf Definition 3.2.}  Let $(H,\alpha)$  be a monoidal Hom-Hopf algebra and  $(A,\b)$  a right
$(H,\a)$-Hom-comodule algebra.  A relative Hom-Hopf module $(M,\mu)$
is a right $(A,\b)$-Hom-module which is also a right
$(H,\a)$-Hom-comodule with the  coaction structure $ \rho _{M} : M
\rightarrow M \otimes  H $ defined by $\rho _{M}(m)=m_{[0]}\o
m_{[1]}$ such that the following compatible condition holds: for all
$m\in M $ and $a\in A$,
\begin{eqnarray*}
&&\rho _{M}( m a) = m _{[0]} \c a _{[0]} \otimes  m _{[1]} a _{[1]}.
\end{eqnarray*}

A morphism between two right relative Hom-Hopf modules is a $k$-linear map
which is a morphism in the categories $\widetilde{\mathscr{H}}(\mathscr{M}_k)_A$ and $\widetilde{\mathscr{H}}(\mathscr{M}_k)^H$ at the same time.  $ \widetilde{\mathscr{H}}(\mathscr{M}_k)^{H}_{A}$ will denote  the category of right relative Hom-Hopf modules and morphisms
between them. \\

\noindent{\bf Proposition 3.3.}
 The forgetful functor
$F:  \widetilde{\mathscr{H}}(\mathscr{M}_k)^{H}_{A} \rightarrow \widetilde{\mathscr{H}}(\mathscr{M}_k)_{A}$ has a right
 adjoint $G:  \widetilde{\mathscr{H}}(\mathscr{M}_k)_{A}\rightarrow \widetilde{\mathscr{H}}(\mathscr{M}_k)^{H}_{A}$.
 $G$ is defined by
$$
 G(M) = M \otimes H,
 $$
 with structure maps
\begin{eqnarray*}
( m  \otimes h  ) \cdot a =   m \c a _{[0]}\otimes h  a _{[1]},
\end{eqnarray*}
\begin{eqnarray*}
\rho _{G(M)} ( m  \otimes h ) =  (\mu^{-1}( m ) \otimes h _{(1)}) \otimes \alpha(h _{(2)}),
\end{eqnarray*}
for all $a \in A $ and $m\in M, h\in H$.

{\it Proof.}
Let us first show that $G (M)$ is an object of $ \widetilde{\mathscr{H}}(\mathscr{M}_k)^{H}_{A}$ . It is routine to check that $G (M)$ is a  right $(H,\a)$-Hom-comodule and a right $(A,\b)$-Hom-module. Now we only check the compatibility condition, for all $a \in A $.
 Indeed,
\begin{eqnarray*}
 &  & \rho_{G(M)}  ((m  \otimes h)\cdot a)\\
 & =&  \rho_{G(M)}  ( m \c a _{[0]} \otimes h a _{[1]})\\
 & = & \mu^{-1}(m)\c \b^{-1}(a_{[0]})\o h_{(1)}a _{[1](1)}\o \alpha(h_{(2)}a _{[1](2)})\\
 & = & \mu^{-1}(m)\c a_{[0][0]}\o h_{(1)}a _{[0][1]}\o \alpha(h_{(2)})a _{[1]}\\
 & = & (m  \otimes h)_{[0]}\c a_{[0]}\o (m  \otimes h)_{(1)}a _{[1]}\\
 &=& \rho (m\otimes c )\cdot a.
  \end{eqnarray*}
This is exactly what we have to show.

For an $A$-linear map $\varphi : (M, \mu) \rightarrow (N, \nu)$, we put
 $$
 G (\varphi)  = \varphi \otimes id _{H}:
M \otimes H\rightarrow N\otimes H .
 $$
Standard computations show that $G (\varphi)$ is  morphisms of right $(A,\b)$-Hom-modules and right $(H,\a)$-Hom-comodules. Let us describe the
 unit $\eta$ and the counit $\delta$ of the adjunction. The unit is described by the coaction:
 for $M \in   \widetilde{\mathscr{H}}(\mathscr{M}_k)^{H}_{A} $,  we define
 $\eta _{M} :
  M \rightarrow M \otimes H$  as follows:
 for all $m \in M $,
 $$
 \eta _{M}(m) =  m _{[0]}
 \otimes m _{[1]}.
 $$
 We can check that $\eta_{M} \in  \widetilde{\mathscr{H}}(\mathscr{M}_k)^{H}_{A}$ . For any $N \in \widetilde{\mathscr{H}}(\mathscr{M}_k)_{A}$, we define
 $\delta _{N}: N\otimes H\rightarrow N$, for all $n \in N$ and $h \in H $,
 $$
 \delta _{N}(n  \otimes h  )
 = \varepsilon(h )n ,
 $$
 we can check that $\delta _{N}$ is $(A,\b)$-linear.
It is easy to check that $\eta _{M}\in   \widetilde{\mathscr{H}}(\mathscr{M}_k)^{H}_{A}$ .
We can check that $\eta$ and $\delta$ defined above are all natural transformations and
 satisfied
 $$
 G (\delta _{N})\circ \eta _{G (N)} = I _{G (N)},
 $$
 $$
 \delta _{F (M)}\circ F (\eta _{M}) = I _{F (M)},
 $$
 for all $M \in  \widetilde{\mathscr{H}}(\mathscr{M}_k)^{H}_{A}$ and $N \in  \widetilde{\mathscr{H}}(\mathscr{M}_k)_{A}$ .

\section{ Structure type theorem and injective type properties  for relative Hom-Hopf modules}
\noindent{\bf Definition 4.1.} Let $(H, \a)$ be a Monoidal Hom Hopf algebra and  $(A,\b)$  a right
$(H,\a)$-Hom-comodule algebra. The map $\phi: (H, \a)\rightarrow (A,\b)$ is called a  total integral such that the following conditions
are satisfied:
$$\rho_A\phi=(\phi\o id_H)\Delta_H,~~~ \phi\a=\b\phi,~~~~~\phi(1_H)=1_A.$$

\noindent{\bf Lemma 4.2.} Let $(H, \a)$ be a monoidal Hom-Hopf algebra and  $(A,\b)$  a right
$(H,\a)$-Hom-comodule algebra with a total integral $\phi: (H, \a)\rightarrow (A,\b)$ and $M\in
\widetilde{\mathscr{H}}(\mathscr{M}_k)^{H}_{A}$,
\begin{eqnarray*}
\lambda_M: M\o H\rightarrow M, ~~~~~~~~~~~~~~ m\o h\mapsto \mu^{-1}(m_{[0]}) \c \phi(S(m_{[1]})\a(h)).
\end{eqnarray*}
Then the following assertions hold:

(1) $\lambda_M\rho_M=id_M$,

(2) $\lambda_M$ is  a morphism of right $(H,\a)$-Hom-comodules, and the right $(H,\a)$-Hom-coaction on $M\o H$ given by $\rho(m\o h)=(\mu(m)\o h_{(1)})\o \a^{-1}(h_{(2)})$ for any $m\in M$ and $h\in H$,

(3)if $\phi: (H, \a)\rightarrow (Z(A), \b)$ (the center of $A$) is a multiplication map, then $\lambda_M$ is a morphism
in $ \widetilde{\mathscr{H}}(\mathscr{M}_k)^{H}_{A}$ .

{\it Proof. }(1) For any $m\in M$, we have
\begin{eqnarray*}
\lambda_M\rho_M(m)&=&\lambda_M(m_{[0]}\o m_{[1]})\\
&=& \mu^{-1}(m_{[0][0]}) \c \phi(S(m_{[0][1]})\a(m_{[1]}))\\
&=& m_{[0]} \c \phi(S(m_{[1](1)})m_{[1](2)})\\
&=&  m_{[0]} \c \phi(\varepsilon(m_{[1]}))\\
&=& \mu^{-1}(m) \c 1_A=m.
\end{eqnarray*}
(2) For any $m\in M$ and $h\in H$, we have
\begin{eqnarray*}
&&\rho_M\lambda_M(m\o h)\\
&=&\rho_M(\mu^{-1}(m_{[0]}) \c \phi(S(m_{[1]})\a(h)))\\
&=& \mu^{-1}(m_{[0][0]}) \c \phi(S(m_{[1](2)})\a(h_{(1)}))\o  \a^{-1}(m_{[0][1]})(S(m_{[1](1)})\a(h_{(2)}))\\
&=& \mu^{-2}(m_{[0]}) \c \phi(\a(S(m_{[1](2)(2)}))\a(h_{(1)}))\o  \a^{-1}(m_{[1](1)}) (\a(S(m_{[1](2)(1)}))\a(h_{(2)})\\
&=& \mu^{-2}(m_{[0]}) \c \phi(S(m_{[1](2)})\a(h_{(1)}))\o m_{[1](1)(1)} (\a(S(m_{[1](1)(2)}))\a(h_{(2)}))\\
&=& \mu^{-2}(m_{[0]}) \c \phi(S(m_{[1](2)})\a(h_{(1)}))\o ( \a(m_{[1](1)(1)})
\a(S(m_{[1](1)(2)})))h_{(2)}\\
&=& \mu^{-2}(m_{[0]}) \c \phi(\a^{-1}(S(m_{[1]}))\a(h_{(1)}))\o \a(h_{(2)})\\
&=& (\lambda_M\o id_H)((\mu^{-1}(m)\o h_{(1)})\o \a(h_{(2)}))\\
&=& (\lambda_M\o id_H)\rho_{M\o H}(m\o h).
\end{eqnarray*}
(3) For any $m\in M,~h\in H$ and $b\in A$, we have
\begin{eqnarray*}
&&\lambda_M((m\o h)\c b)\\
&=& \lambda_M(m\c b_{[0]}\o hb_{(1)})\\
&=& \mu^{-1}(m_{[0]}\c b_{[0][0]}) \c \phi(S(m_{[1]}b_{[0][1]})\a(h b_{[1]}))\\
&=& \mu^{-1}(m_{[0]}\c b_{[0][0]}) \c \phi(S(m_{[1]})S(b_{[0][1]})\a(h b_{[1]}))\\
&=& \mu^{-1}(m_{[0]}\c b_{[0][0]}) \c \phi(\a(S(m_{[1]})[S(b_{[0][1]})h b_{[1]}])\\
&=& \mu^{-1}(m_{[0]}\c b_{[0][0]}) \c \phi(\a(S(m_{[1]})[S(b_{[0][1]}) (b_{[1]}h)])\\
&=& \mu^{-1}(m_{[0]}\c b_{[0][0]}) \c \phi(\a(S(m_{[1]})[(\a^{-1}(S(b_{[0][1]}))b_{[1]})\a(h)])\\
&=& (\mu^{-1}(m_{[0]})\c b_{[0]}) \c \phi(\a(S(m_{[1]})[(\a^{-1}(S(b_{[1](1)}))\a^{-1}(b_{[1](2)}))\a(h)])\\
&=& (\mu^{-1}(m_{[0]})\c \b^{-1}(b)) \c \phi(\a(S(m_{[1]})\a^{2}(h))\\
&=& m_{[0]}\c (\b^{-1}(b) \phi(S(m_{[1]})\a(h)))\\
&=& m_{[0]}\c (\phi(\a^{-4}(S(m_{[1]}))\a^{-3}(h))\b^{-1}(b))\\
&=& (\mu^{-1}(m_{[0]})\c \phi(S(m_{[1]})\a(h))) \c b\\
&=& \lambda_M(m\o h)\c b.
\end{eqnarray*}

\noindent{\bf Proposition  4.3.} Let $(H, \a)$ be a monoidal Hom-Hopf algebra and  $(A,\b)$  a right $(H,\a)$-Hom comodule
algebra with a total integral $\phi: (H, \alpha)\rightarrow (A, \beta)$. Then $M\in  \widetilde{\mathscr{H}}(\mathscr{M}_k)^{H}_{A} $ is injective as a right $(H, \a)$-Hom-comodule.\\

If $H$ is a Hopf algebra, then we obtain the main result of [3, Theorem 1].\\

\noindent{\bf Corollary 4.4.} Let $H$ be a Hopf algebra and $A$ a
right $H$-comodule algebra. If there is a right $H$-comodule map
$\phi: (H, \alpha)\rightarrow (A, \beta)$ such that $\phi(1_H)=1_A$,
then every relative $(H,A)$-Hopf-module
is injective as a right $H$-comodule.\\

Let $M$ be a relative Hom-Hopf module, and let
\begin{eqnarray*}
M_{0}=\{m\in M \mid \rho_{M}(m)=\mu^{-1}(m)\o 1_H\}
\end{eqnarray*}
be an invariant subspace of $M$ and $M_0$ is a right $(C, \b)$-Hom-module, where
\begin{eqnarray*}
C=\{b\in A\mid \rho_{A}(b)=\b^{-1}(b)\o 1_H\}
\end{eqnarray*}
is a subalgebra of $A$.\\

\noindent{\bf Proposition  4.5.} Let $(H, \a)$ be a monoidal Hom-Hopf algebra and  $(A,\b)$  a right
$(H,\a)$-Hom-comodule algebra with a total integral $\phi: (H, \a)\rightarrow (A,\b)$ and $M\in
\widetilde{\mathscr{H}}(\mathscr{M}_k)^{H}_{A}$. Assume that $\phi$
is a multiplication map and let
$$\tau_{M}: (M, \mu)\rightarrow (M, \mu)$$
be the trace map defined by
$$m\mapsto m_{[0]}\c \phi(S(m_{[1]})).$$
Then the following assertions hold:

(1) $\tau_{M}(m)\in M_{0}$ and $\tau\mid_{M_0}=id$,

(2)$\tau_A: (A, \b)\rightarrow (C, \b)$ defined by $b\mapsto b_{[0]}\phi(S(b_{[1]}))$ is a morphism of left $(C,\b)$-Hom-modules, so that $(C,\b)$ is a direct summand of $(A,\b)$ as a sum of left $(C,\b)$-Hom-modules,

(3) if $Im\phi\subseteq Z(A)$, then $\tau_{M}: (M, \mu)\rightarrow (M, \mu)$ is a morphism of right $(C,\b)$-Hom-modules.

The exact sequence

$$(M, \mu)\stackrel{\tau_{M}}\longrightarrow (M_0, \mu)\longrightarrow 0,$$
thus obtained splits as a sequence of right $(C,\b)$-Hom-modules.

{\it Proof. }(1) For any $m\in M$, we have
\begin{eqnarray*}
\rho(\tau_{M}(m))&=&\rho(m_{[0]} \phi(S(m_{[1]})))\\
&=& m_{[0][0]} \phi(S(m_{[1](2)})) \o m_{[0][1]} \phi(S(m_{[1](1)}))\\
&=&\mu^{-1}(m_{[0]}) \phi(\a(S(m_{[1](2)(2)}))) \o m_{[1](1)}\phi(\a(S(m_{[1](2)(1)})))\\
&=&\mu^{-1}(m_{[0]}) \phi(S(m_{[1](2)})) \o \a(m_{[1](1)(1)}) \phi(\a(S(m_{[1](1)(2)})))\\
&=& \mu^{-1}(m_{[0]}) \phi(\a^{-1}(S(m_{[1]}))) \o 1_H\\
&=& \mu^{-1}(\tau_{M}(m))\o 1_H.
\end{eqnarray*}

For any $n\in M_0$,
\begin{eqnarray*}
\tau_{M}(n)&=& n_{[0]}\phi(S(n_{[1]}))\\
&=& \mu^{-1}(n)1_A=n.
\end{eqnarray*}

(2) For any $c\in C$ and $a\in A$,
\begin{eqnarray*}
 \tau_A(ca)&=& (c_{[0]}a_{[0]})\phi(S(c_{[1]}a_{[1]}))\\
   &=& (\b^{-1}(c)a_{[0]}) \phi(\a(S(a_{[1]})))\\
   &=& c(a_{[0]}\c \phi(S(a_{[1]})))=c\tau_A(a),
 \end{eqnarray*}
thus, $\tau_A: (A, \b)\rightarrow (C, \b)$  is a morphism of  left $(C,\b)$-Hom-modules, and by (1), $(C, \b)$ is a direct summand of $(A,\b)$ as a sum of left $(C,\b)$-Hom-modules.

(3)  For any $c\in C$ and $m\in M$,
\begin{eqnarray*}
\tau_{M}(m\c c) &=& (m_{[0]}\c c_{[0]})\phi(S(m_{[1]}c_{[1]}))\\
&=& (m_{[0]}\c \b^{-1}(c))\phi(\a(S(m_{[1]})))\\
&=& \mu(m_{[0]})\c (\b^{-1}(c))\phi(S(m_{[1]})))\\
&=& \mu(m_{[0]})\c (\phi(S(m_{[1]}))\b^{-1}(c)))\\
&=& (m_{[0]}\c \phi(S(m_{[1]})))\c c=\tau_{M}(m)\c c.
\end{eqnarray*}
Thus, $\tau_{M}$ is  a morphism of right  $(C,\b)$-Hom-modules, and by (1),  the exact sequence

$$(M, \mu)\stackrel{\tau_{M}}\longrightarrow (M_0, \mu)\longrightarrow 0.$$
Thus obtained splits as a sequence of right $(C,\b)$-Hom-modules.

 \section{\textbf{A Maschke-type theorem for relative Hom-Hopf modules}}
\def\theequation{5. \arabic{equation}}
\setcounter{equation} {0} \hskip\parindent

In this section, we shall  give necessary and sufficient conditions for the functor $F$
 which forget the $(H, \a)$-coaction to be separable,  we prove a Maschke type theorem  for relative Hom-Hopf modules as an application.\\

\noindent{\bf Definition 5.1.} Let $(H, \a)$ be a monoidal Hom-Hopf algebra and  $(A,\b)$  a right
$(H,\a)$-Hom-comodule algebra. A $k$-linear map
  $$
\theta : (H, \a)  \otimes (H, \a)  \rightarrow  (A, \b)
  $$ such that $\theta\circ (\a\o \a)=\b\circ \theta$  is called a  normalized $(A, \b)$-integral, if $\theta$ satisfies the following conditions:

  (1) For all $h,g \in H $,
 \begin{eqnarray}
 && \theta (\a^{-1}(g) \otimes h _{(1)} )\otimes \a(h _{(2)}) \nonumber\\
 &=& \beta( \theta  (g _{(2)} \otimes \a^{-1}(h))
   _{[0]})\otimes g _{(1)} \theta (g _{(2)} \otimes \a^{-1}(h))_{[1]}.
 \end{eqnarray}

(2) For all  $h \in H $,
 \begin{equation}
 \theta (h _{(1)} \otimes h _{(2)})
 = 1 _{A } \varepsilon (h).
 \end{equation}

(3) For all $a \in A, h,g \in H$,
 \begin{equation}
  \b^2 (a _{[0][0]})
   \theta (\a^{-1}(g) a _{[0][1]} \otimes \a^{-1}(h)\a^{-1}(a _{[1]}))
   = \theta ( g \otimes h ) a.
 \end{equation}

\noindent{\bf Theorem 5.2.}
Let $(H, \a)$ be a monoidal Hom-Hopf algebra and  $(A,\b)$  a right $(H,\a)$-Hom-comodule algebra, the following assertions are equivalent,

(1) The left adjoint $F$ in Proposition 3.3 is separable,

(2) There exists a normalized $(A, \b)$-integral
$\theta : (H, \a)  \otimes (H, \a)  \rightarrow  (A, \b) $.

{\it Proof.} $(2) \Longrightarrow(1)$. For any relative Hom-Hopf module $M$, we define
\begin{eqnarray*}
&  & \nu_{M} : M  \otimes H \rightarrow M ,\\
&  &\nu_{M} (m \otimes h  )\\
&=& \mu(m _{[0]}) \theta  (m _{[1]} \otimes \a^{-1}(h)),
\end{eqnarray*}
for all  $m \in M $ and $h \in H $.
Now, we shall check that $\nu _{M} \in \widetilde{\mathscr{H}}(\mathscr{M}_k)^{H}_{A}$. In fact, for all $m \in M$, $h \in H$ and $a \in A $, it is easy to get that
$$
\nu_M(\mu(m) \o \a(h)) = \mu(\nu_M(m \o h)).
$$
We also have
\begin{eqnarray*}
&&\nu _{M} ((m  \otimes h )\cdot a)\\
&=&\nu_{M}( m a _{[0]} \otimes h  a _{[1]})\\
&=& (\mu(m_{[0]})\cdot\b (a _{[0][0]})) \theta
    (m_{[1]} a _{[0][1]}\otimes \a^{-1}(h)  \a^{-1}(a _{[1]}))\\
&=& \mu^2(m_{[0]})\cdot(\b(a _{[0][0]})\b^{-1}(\theta
    (m_{[1]} a _{[0][1]}\otimes \a^{-1}(h)  \a^{-1}(a _{[1]})))\\
&=& \mu^2(m_{[0]})\cdot(\b(a _{[0][0]})\theta(\a^{-1}(m_{[1]})\a^{-1}(a _{[0][1]}) \otimes \a^{-2}(h)  \a^{-2}(a _{[1]})))\\
&\stackrel{(5.3)}{=}& \mu^2(m_{[0]})\cdot(\theta(m_{[1]} \otimes \a^{-1}(h)) \b^{-1}(a))\\
&=& (\mu(m_{[0]})\cdot \theta(m_{[1]} \otimes \a^{-1}(h)))\cdot a\\
&=& (\nu _{M}(m  \otimes h ))\cdot a.
\end{eqnarray*}
Hence  it is a morphism of $(A, \b)$-Hom-modules. Next, we shall check that $\nu _{M}$ is a morphism of Hom-comodules
over $(H, \a)$. It is sufficient to check that
 $$
 \rho _{M}  \circ \nu_{M}  = (\nu_{M}  \otimes id _{H})
 \circ \rho _{M}
 $$
 holds. For all $m \in M $ and $h \in H$,  we have
\begin{eqnarray*}
&&\rho _{M}\circ \nu_{M}(m  \otimes h )\\
&=& \rho _{M}  (\mu(m_{[0]} )
\theta (m_{[1]} \otimes \a^{-1}(h )  ))\\
&=& (\mu(m_{[0]})
\theta(m_{(1)} \otimes \a^{-1}(h))_{[0]}
 \otimes (\mu(m_{[0]}) \theta (m_{[1]} \otimes \a^{-1}(h) )) _{[1]} \\
&=& \mu(m_{[0][0]}) \theta (m_{[1]}    \otimes \a^{-1}(h))_{[0]}   \otimes
   \a( m_{[0][1]} )
\theta (m_{(1)}\otimes \a^{-1}(h) ) _{[1]} \\
&=& m_{[0]}
\theta (\a(m_{[1](2)})\otimes \a^{-1}(h))_{[0]}
 \otimes\a(m_{[1](1)}) \theta (\a(m_{[1](2)})\otimes \a^{-1}(h)  ) _{[1]}\\
&\stackrel{(5.1)}{=}& m_{[0]}
\b^{-1}(\theta (m_{[1]} \otimes h_{(1)}))
 \otimes \a(h  _{(2)}) \\
& = & m_{[0]}
\theta (\a^{-1}(m_{[1]}) \otimes \a^{-1}(h_{(1)}))
 \otimes \a(h  _{(2)}) \\
& = &(\nu_{M}  \otimes id _{H})
 \circ \rho _{M} (m \otimes h).
\end{eqnarray*}
For all $m \in M $, since
\begin{eqnarray*}
& & \nu _{M}  \circ \eta _{M}  (m) = \nu _{M} (m _{[0]}
 \otimes m _{[1]})\\
&= & \mu(m _{[0][0]}) \theta  (m _{[0][1]}
 \otimes  \a^{-1}(m _{[1]}))\\
&= &  m _{[0]} \theta (m _{[1](1)}\otimes m_{[1](2)})\stackrel{(5.2)}{=}
 m.
\end{eqnarray*}
So the left adjoint $F$ in Proposition 3.3 is separable follows by Rafael theorem.

$(1) \Longrightarrow (2)$. We consider the following  relative Hom-Hopf module $A\o H $,
and the $(A, \b)$-actions and $(H, \a)$-coaction are defined as follows:
$$
\left\{
  \begin{array}{ll}
(a \otimes h ) \cdot b = a b _{[0]}\otimes h b _{[1]};\\
  \rho_{A\o H}  ( a \otimes h  )
  =(\b^{-1}(a) \otimes h _{(1)})\otimes \a(h _{(2)}),
\end{array}
\right.
$$
for any $a, b\in A$ and $h\in H$.

Evaluating at this  object, the retraction $\nu$ of the unit of the adjunction in Proposition 3.3 yields a morphism
$$
\nu _{A\o H  }:  (A  \otimes H) \otimes H  \rightarrow A \otimes H
$$
such that, for all  $a \in A , h \in H$,

$$
\nu _{A \o H}((a \otimes h _{(1)})
\otimes h _{(2)})= a \otimes h.$$

 It can be used to construct $\theta$ as follows:
$$
\theta :H \otimes H\rightarrow A,
$$
$$
\theta (h \otimes g) = r_A (id _{A }\otimes \varepsilon) \nu _{A \o H}
((1 _{A }
\otimes h) \otimes g ),
$$
where $r$ means the right unit constraint. For all $h \in H $, since
\begin{eqnarray*}
&&\theta (h _{(1)}\otimes h _{(2)})\\
&=&r_A(id _{A }\otimes \varepsilon) \nu _{A \o H}((1 _{A}
\otimes h _{(1)})
\otimes h _{(2)})\\
&=&r_A(id _{A }\otimes \varepsilon) (1 _{A}\otimes h )
= 1 _{A} \varepsilon(h).
\end{eqnarray*}
Hence condition Eq.(5.2) follows. It can be seen to obey Eq.(5.3) by naturality and the
$(A, \b)$-modules map of $\nu$.

 The verification of Eq.(5.1) is more
involved. For any right $(H, \a)$-Hom-comodule $M$, we consider the  relative Hom-Hopf module
$M \otimes  A$,
the $(A, \b)$-action and $(H,\a)$-coaction are defined as follows: for all $m \in M$ and $a, b \in A $,
$$
\left\{
  \begin{array}{ll}
    (m \otimes a) \cdot b = \mu^{-1}(m) \otimes a\b(b), \\
  \rho_{M\o A} ( m \otimes a)= (m _{[0]} \otimes a _{[0]})\otimes m _{[1]}  a _{[1]}.
  \end{array}
\right.
$$
In particular, there is a  relative Hom-Hopf module $H \otimes  A$ and the map
$$
\xi : H \otimes A  \rightarrow A  \otimes H\
$$
$$
\xi(h \otimes a) = \b(a _{[0]}) \otimes \a^{-1}(h) a _{[1]}.
$$
Since $\xi$ is both right $(A, \b)$-linear and right $(H, \a)$-colinear, thus we have
\begin{eqnarray}
\xi(\nu_{H \o A}((h \o a) \o g)) &=& \nu_{A \o H}((\xi \o id_H)((h \o a) \o g)) \nonumber\\
&=& \nu_{A \o H}((\b(a_{[0]}) \o \a^{-1}(h)a_{[1]})\o g).
\end{eqnarray}
It is not hard to check that $GF(H \o A)= (H\o A) \o H \in {}^{H}\widetilde{\mathscr{H}}(\mathscr{M}_k)^{H}_{A}$, and its left $(H, \a)$-Hom comodule structure
is given by
$$
(h \o a) \o g \mapsto \a(h_{(1)}) \o ((h_{(2)} \o \b^{-1}(a)) \o \a^{-1}(g)).
$$
Also $H \o A\in {}^{H}\widetilde{\mathscr{H}}(\mathscr{M}_k)^{H}_{A}$, and the left $(H,\a)$-coaction of $H \o A$ is given by
$$
h \o a\mapsto \a(h_{(1)}) \o (h_{(2)} \o \b^{-1}(a)).
$$
We also get $\nu_{H\o A}:(H\o A) \o H\rightarrow H \o A$ is a Hom morphism in ${}^{H}\widetilde{\mathscr{H}}(\mathscr{M}_k)^{H}_{A}$, which means
\begin{eqnarray}
&\nu_{H\o A}((h \o a) \o g)_{[-1]} \o \nu_{H\o A}((h \o a) \o g)_{[0]} \nonumber\\
&~~~~~~~~= \a(h_{(1)}) \o \nu_{H\o A}((h_{(2)} \o \b^{-1}(a)) \o \a^{-1}(g)).
\end{eqnarray}
Thus we conclude that $\nu_{H\o A}$ is left and right $(H,\a)$-colinear. Take $h, g\in H$, and put
$$\nu_{A\o H}((1_A \o h) \o g) =\sum_{i} a_i\o q_i \in A \o H,$$
$$\nu_{H\o A}((h\o 1_A)\o g)=\sum_{i} p_i\o b_i \in H \o A,$$
we obtain
\begin{eqnarray*}
&&\b(\theta(h_{(2)} \o \a^{-1}(g))_{[0]}) \o h_{(1)}\theta(h_{(2)} \o \a^{-1}(g))_{[1]}\\
&=& \b(r_A(id_A \o \varepsilon)\nu_{A \o H}((1_A \o h_{(2)}) \o \a^{-1}(g))_{[0]}) \o h_{(1)}\\
&~~~~&\cdot(r_A(id_A \o \varepsilon)\nu_{A \o H}((1_A \o h_{(2)}) \o \a^{-1}(g)))_{[1]} \\
&\stackrel{(5.4)}{=}& \b(r_A(id_A \o \varepsilon)\xi\nu_{H \o A}(( h_{(2)} \o 1_A) \o \a^{-1}(g))_{[0]}) \o h_{(1)}\\
&~~~~&\cdot(r_A(id_A \o \varepsilon)\xi\nu_{H \o A}(( h_{(2)} \o 1_A) \o \a^{-1}(g))_{[1]}) \\
&\stackrel{(5.5)}{=}& \sum_{i} \b({r_A(id_A \o \varepsilon)\xi({p_i}_{(2)} \o \b^{-1}(b_i))}_{[0]}) \o {p_i}_{(1)}({r_A(id_A \o \varepsilon)\xi({p_i}_{(2)} \o \b^{-1}(b_i))}_{[1]})\\
&=&\sum_{i} \b(r_A(id_A \o \varepsilon)({b_i}_{[0]} \o \a^{-1}({p_i}_{(2)}){b_i}_{[1]})_{[0]}) \o {p_i}_{(1)}(r_A(id_A \o \varepsilon)({b_i}_{[0]} \o \a^{-1}({p_i}_{(2)}){b_i}_{[1]})_{[1]})\\
&=& \sum_{i}\b({b_i}_{[0]}) \o {p_i}_{(1)}\varepsilon({p_i}_{(2)})({b_i}_{[1]})\\
&=& \sum_{i}\xi(p_i \o b_i)
= \xi(\nu_{H \o A}((h \o 1_A) \o g)).
\end{eqnarray*}
Use the fact that $\nu_{A \o H}$ is a morphism of right $(H, \a)$-Hom comodules, we also have
\begin{eqnarray*}
&&\theta(\a^{-1}(h) \o g_{(1)}) \o \a(g_{(2)})\\
&=& r_A(id_A \o \varepsilon)\nu_{A \o H}((1_A \o \a^{-1}(h))\o g_{(1)})\o \a(g_{(2)})\\
&=& \sum_{i}r_A(id_A \o \varepsilon)(\b^{-1}(a_i) \o {q_i}_{(1)}) \o \a({q_i}_{(2)})\\
&=& \sum_{i} a_i \o q_i = \nu_{A \o H}((1_A \o h) \o g)\\
&\stackrel{(5.4)}{=}& \xi(\nu_{H \o A}((h \o 1_A) \o g)).
\end{eqnarray*}
Hence, we can get condition Eq.(5.1).

We will now investigate the relation between the total integrals and
the normalized $(A,\b)$-integrals. This will explain our terminology, and we will also
prove that the forgetful functor is separable if and only if there exists a total
integral $\phi: (H, \alpha)\rightarrow (A, \beta)$ such that the image of $\rho_A \circ\phi$ is contained in the center of $H\o A$.\\

 \noindent{\bf Proposition 5.3.} Let $(H, \a)$ be a monoidal Hom-Hopf algebra and $(A, \b)$ a  right $(H, \a)$-Hom-comodule algebra. If
$\theta: (H, \a) \otimes (H, \a) \rightarrow (A, \b)$ is a normalized $(A,\b)$-integral for $(H, A,H)$, then the $k$-linear map
$$
\phi: (H, \a)\rightarrow (A, \b),\ \phi(h) = \theta(1 _{H}\otimes h),
$$
for all $h \in H$,  is a  total integral.

{\it Proof.} Notice first that
$\phi(1_{H}) = \theta(1 _{H}\otimes 1 _{H}) =\varepsilon _{H}(1_{H})1_{A}=1_{A}$. Since
\begin{eqnarray*}
 && \theta (\a^{-1}(g) \otimes \a^{-1}(h _{(1)}) )\otimes \a(h _{(2)})\\
 &=& ( \theta  (\a(g _{(2)}) \otimes \a^{-1}(h)))
   _{(0)}\otimes\a(g _{(1)})( \theta (\a(g _{(2)}) \otimes \a^{-1}(h)))_{(1)}.
\end{eqnarray*}
It follows by taking $g=1_{H}$ that
$$
\theta  (1_{H} \otimes \a^{-1}(h_{1 }))\otimes \a(h _{2})=\theta  (1_{H}  \otimes \a^{-1}(h))_{[0]} \otimes \a(\theta  (1_{H} \otimes \a^{-1}(h))_{[1]}),
$$
and applying $\a \o \a^{-1}$ to the above identity, we have
$$
\theta  (1_{H} \otimes \a^{-1}(h_{1 }))\otimes h _{2}=\theta  (1_{H}  \otimes \a^{-1}(h))_{[0]} \otimes \theta  (1_{H} \otimes \a^{-1}(h))_{[1]}.
$$
So we obtain
$$
\phi(h_{1 })\otimes h _{2}=\phi(h)_{[0]} \otimes \phi  (h)_{[1]}.
$$
It is easy to check that $\phi\a=\b\phi$. So $\phi$ is a  total integral.

Let $\phi: (H, \a)\rightarrow (A, \b)$ be  a total integral  for the right $(H, \a)$-Hom-comodule algebra $(A, \b)$, and define
$$
\theta: (H, \a) \otimes (H, \a) \rightarrow (A, \b),~~~~~~~ \theta(g\otimes h)= \phi(hS^{-1}(g)),
$$
for all $g, h \in H$.\\

\noindent{\bf Theorem 5.4.}
Let $(H, \a)$ be a monoidal Hom-Hopf algebra and $(A, \b)$ a  right $(H, \a)$-Hom-comodule algebra, and $\phi: (H, \a)\rightarrow (A, \b)$  a  total integral.
If
$$
g \phi(h)_{[1]} \otimes \phi(h)_{[0]} = \phi(h)_{[1]}g \otimes \phi(h)_{[0]},~~~\phi(h) \in Z(A).
$$
Then $\theta$ is a normalized $(A, \b)$-integral.

{\it Proof.} For any $h, g\in H$ and $a\in A$, we have
\begin{eqnarray*}
&&\b^2 (a _{[0][0]})
   \theta (\a^{-1}(g) a _{[0][1]} \otimes \a^{-1}(h)\a^{-1}(a _{[1]}))\\
   &=& \b (a _{[0]})
   \theta (\a^{-1}(g) a _{[1](1)} \otimes \a^{-1}(h)a _{[1](2)})\\
   &=& \b (a _{[0]})\phi(\a^{-1}(h)a _{[1](2)}S^{-1}(\a^{-1}(g) a _{[1](1)}))\\
   &=&\b (a _{[0]})\phi(h[(\a^{-1}(a _{[1](2)}S^{-1}(\a^{-1}(a _{[1](1)}))))S^{-1}(\a^{-1}(g)) ])\\
   &=& a\phi(hS^{-1}(g))\\
   &=& \theta ( g \otimes h ) a,
\end{eqnarray*}

and
\begin{eqnarray*}
 && \beta( \theta  (g _{(2)} \otimes \a^{-1}(h))
   _{[0]})\otimes g _{(1)} \theta (g _{(2)} \otimes \a^{-1}(h))_{[1]}\\
 &=& \beta( \phi (\a^{-1}(h)S^{-1}(g _{(2)}))
   _{[0]})\otimes  \phi (\a^{-1}(h)S^{-1}(g _{(2)}))_{[1]}g _{(1)}\\
   &=& \phi (h_{(1)}S^{-1}(\a(g _{(2)(2)})))\otimes  (\a^{-1}(h_{(2)})S^{-1}(g _{(2)(1)}))g _{(1)}\\
   &=& \phi (h_{(1)}S^{-1}(g _{(2)}))\otimes  (\a^{-1}(h_{(2)})S^{-1}(g _{(1)(2)}))\a(g _{(1)(1)})\\
    &=& \phi (h_{(1)}S^{-1}(g _{(2)}))\otimes  h_{(2)}(S^{-1}(g _{(1)(2)})g _{(1)(1)})\\
    &=&\phi (h_{(1)}S^{-1}(\a^{-1}(g)))\otimes  \a(h_{(2)})\\
 &=& \theta (\a^{-1}(g) \otimes h _{(1)} )\otimes \a(h _{(2)}),
 \end{eqnarray*}
\begin{eqnarray*}
\theta(h _{1} \otimes h _{2}) =  \varphi(h _{2}S ^{-1}(h _{1}))=\varepsilon _{H}(h)1_{A}.
\end{eqnarray*}
It is easy to check that $\phi\a=\b\phi$. So $\theta$ is a normalized $(A, \b)$-integral.

Since separable functors reflect well the semisimplicity of the objects of a categogy , by Theroem 5.2,  we will prove a Maschke type theorem for   relative Hom-Hopf modules as an application.\\

\noindent{\bf Lemma 5.5.} Let $(H, \a)$ be a monoidal Hom-Hopf algebra and  $(A,\b)$  a right
$(H,\a)$-Hom-comodule algebra with a total integral $\phi: (H, \a)\rightarrow (A, \b)$ and $(M,
\mu),  (N, \nu) \in \widetilde{\mathscr{H}}(\mathscr{M}_k)^{H}_{A}$
and a Hom-morphism $f: (N, \nu)\rightarrow (M, \mu) $. Let $$f_\phi:
N\stackrel{\rho_N}{\longrightarrow} N\o H\stackrel{f\o
id_H}{\longrightarrow} M\o H\stackrel{\tau}\longrightarrow M,$$ that
is,
$$f_{\phi}(n)=\mu^{-1}(f(n_{[0]})_{[0]})\c \phi(S(f(n_{[0]})_{[1]})\a(n_{[1]})),$$
for any $n\in N$. Then the following assertions hold:

(1) $f_\phi$   is a morphism of right  $(H, \a)$-Hom-comodules,

(2) if $f: (N, \nu)\rightarrow (M, \mu) $ is a  morphism of right  $(A, \b)$-Hom-modules  and $\phi: (H, \a)\rightarrow (Z(A), \b)$ is a multiplication map, then  $f_\phi$ is a morphism of right $(A,\b)$-Hom-module.

{\it Proof. }(1) For any $n\in N$, we have
\begin{eqnarray*}
\rho_M(f_\phi(n))&=& \rho_M(\mu^{-1}(f(n_{[0]})_{[0]})\c \phi(S(f(n_{[0]})_{[1]})\a(n_{[1]})))\\
&=&\mu^{-1}(f(n_{[0]})_{[0][0]})\c \phi(S(f(n_{[0]})_{[1](2)})\a(n_{[1](1)})) \\
&&\hspace{1.5cm}\o \a^{-1}(f(n_{[0]})_{[0][1]})(S(f(n_{[0]})_{[1](1)})\a(n_{[1](2)}))\\
&=&\mu^{-2}(f(n_{[0]})_{[0]})\c \phi(\a(S(f(n_{[0]})_{[1](2)(2)}))\a(n_{[1](1)})) \\
&&\hspace{1.5cm}\o \a^{-1}(f(n_{[0]})_{[1](1)})(\a(S(f(n_{[0]})_{[1](2)(1)}))\a(n_{[1](2)}))\\
&=&\mu^{-2}(f(n_{[0]})_{[0]})\c \phi(S(f(n_{[0]})_{[1](2)})\a(n_{[1](1)})) \\
&&\hspace{1.5cm}\o f(n_{[0]})_{[1](1)(1)}(\a(S(f(n_{[0]})_{[1](1)(2)}))\a(n_{[1](2)}))\\
&=&\mu^{-2}(f(n_{[0]})_{[0]})\c \phi(S(f(n_{[0]})_{[1](2)})\a(n_{[1](1)})) \\
&&\hspace{1.5cm}\o (\a(f(n_{[0]})_{[1](1)(1)})\a(S(f(n_{[0]})_{[1](1)(2)})))n_{[1](2)}\\
&=&\mu^{-2}(f(n_{[0]})_{[0]})\c \phi(\a^{-1}(S(f(n_{[0]})_{[1]}))\a(n_{[1](1)})) \o \a(n_{[1](2)})\\
&=&\mu^{-1}(f(n_{[0][0]})_{[0]})\c \phi(S(f(n_{[0][0]})_{[1]})\a(n_{[0][1]})) \o n_{[1]}\\
&=& (f_\phi\o id_H)\rho_N(n).
\end{eqnarray*}

(2) For any $n\in N$ and $b\in A$, we have
\begin{eqnarray*}
&&f_\phi(n\c b)\\
&=& \mu^{-1}(f(n_{[0]})_{[0]}\c b_{[0][0]})\c \phi(S(f(n_{[0]})_{[1]}b_{[0][1]})\a(n_{[1]}b_{[1]}))\\
&=& \mu^{-1}(f(n_{[0]})_{[0]}\c b_{[0][0]})\c \phi([S(f(n_{[0]})_{[1]} b_{[0][1]})][\a(b_{[1]})\a(n_{[1]}) ])\\
&=& \mu^{-1}(f(n_{[0]})_{[0]}\c b_{[0][0]})\c \phi(\a(S(f(n_{[0]})_{[1]}) [S(b_{[0][1]})(b_{[1]} n_{[1]})])\\
&=& \mu^{-1}(f(n_{[0]})_{[0]}\c b_{[0][0]})\c \phi(\a(S(f(n_{[0]})_{[1]}) [(\a^{-1}(S(b_{[0][1]}))b_{[1]}) \a(n_{[1]})])\\
&=& (\mu^{-1}(f(n_{[0]})_{[0]})\c  b_{[0]})\c \phi(\a(S(f(n_{[0]})_{[1]}) [(\a^{-1}(S(b_{[1](1)}))\a^{-1}(b_{[1](2)})) \a(n_{[1]})])\\
&=& (\mu^{-1}(f(n_{[0]})_{[0]})\c \b^{-1}( b))\c \phi(\a(S(f(n_{[0]})_{[1]})) \a^{2}(n_{[1]}))\\
&=& f(n_{[0]})_{[0]}\c (\b^{-1}( b)\phi(S(f(n_{[0]})_{[1]}) \a(n_{[1]})))\\
&=& f(n_{[0]})_{[0]}\c ( \phi(S(f(n_{[0]})_{[1]}) \a(n_{[1]}))\b^{-1}( b))\\
&=& (\mu^{-1}(f(n_{[0]})_{[0]})\c  \phi(S(f(n_{[0]})_{[1]}) \a(n_{[1]}))) \c b=f_\phi(n)\c b.
\end{eqnarray*}

\noindent{\bf Theorem  5.6.} Let $(H, \a)$ be a monoidal Hom-Hopf algebra and  $(A,\b)$  a right
$(H,\a)$-Hom-comodule algebra with a total integral $\phi: (H, \a)\rightarrow (A, \b)$. If
$\phi: (H, \a)\rightarrow (Z(A), \b)$ is a multiplication map, then
every short exact sequence of  relative Hom-Hopf modules
$$0\longrightarrow (M, \mu)\stackrel{f}\longrightarrow (N, \nu)\stackrel{g}\longrightarrow (P, \pi)\longrightarrow 0$$
which splits as a sequence of $(A,\b)$-Hom-modules also splits as a sequence of relative
Hom-Hopf modules.

\section*{Acknowledgements}

The authors are grateful to the referee for carefully reading the manuscript
and for many valuable comments which largely improved the article. The work was
 supported  by  the NSF of
Jiangsu Province (No. BK2012736),  Southeast University for Postdoctoral Innovation Funds (No. 3207013601)£¬
 and Jiangsu Planned Projects for Postdoctoral Research Funds(No.1302019c).

\end{document}